\documentclass[preprint,12pt]{article}

\usepackage{amssymb,amsmath,array}

\begin{document}

\title{General bound of overfitting for MLP regression models.}

\author{Rynkiewicz, J.}
\date{}
\maketitle

\begin{abstract}
Multilayer perceptrons (MLP) with one hidden layer have been used for a long time to deal with non-linear regression. However, in some task, MLP's are  too powerful models and a small mean square error (MSE) may be more due to overfitting than to actual modelling. If the noise of the regression model is Gaussian, the overfitting of the model is totally determined by the behavior of the likelihood ratio test statistic (LRTS), however in numerous cases the assumption of normality of the noise is arbitrary if not false.  In this paper, we present an universal bound for the overfitting of such model under weak assumptions, this bound is valid without Gaussian or identifiability  assumptions. The main application of this bound is to give a hint about determining the true architecture of the MLP model when the number of data goes to infinite. As an illustration, we use this theoretical result to propose and compare effective criteria to find the true architecture of an MLP.
\end{abstract}



\section{Introduction}
Feed-forward neural networks are well known and popular tools to deal with non-linear regression models.
We can describe MLP models as a parametric family of regression functions. White \cite{White} reviews statistical properties of  MLP estimation in detail. However he leaves an important question pending i.e. the asymptotic behavior of the estimator when the MLP in use has redundant hidden units. If we assume that the noise is Gaussian it is well known that the Least square Estimator (LSE) and the maximum likelihood estimator (MLE) are equivalent and Amari et al. \cite{Amari} give several examples of the behavior of MLE in  case of redundant hidden units. Moreover, if $n$ is the number of observations, Fukumizu \cite{Fukumizu2003} shows that, for unbounded parameters and Gaussian noise, LRTS can have an order lower bounded by $O(\log(n))$  instead of the classical convergence property to a $\chi^2$ law. In the same spirit, Hagiwara and Fukumizu \cite{Hagiwara} investigate relation between LRTS divergence and weight size in a simple neural networks regression problem. Hence, if parameters of MLP models are not bounded, these papers show that, for Gaussian noise, the overfitting is strong, even if the number of data is large. Note that this is no more the case if the parameters of the MLP are supposed to be a priori bounded.

But, even if Gaussian assumption for the noise is standard, it may be not suitable for some models. This assumption is false, for example, when the range of observations is known to be bounded, since Gaussian variables can be arbitrary large in absolute value, even if the probability of such events is small. Hence, we need a theory which gives evaluation of the overfitting of MLP regression without knowing the density of the noise and which works even if the model is not identifiable.  

In this paper, we prove an inequality bounding the MSE difference between the true model and an over-determined model. This inequality shows that, under suitable assumptions, the asymptotic overfitting of the MSE is upper bounded  by the maximum of the square of a Gaussian process. Moreover, this bound shows that suitably penalized MSE criteria allow to select asymptotically the true model. The paper is organized as follows: In section 2 we state the model, section 3 presents our main inequality and in section 4 we apply this inequality to select the optimal architecture of the MLP model. Finally, in  section 5, a little experiment gives us some insight to apply this theoretical results in real life problems.  

\section{The model}
\label{sec:1}
Let  $x=(x(1),\cdots,x(d))^T\in{\mathbb R}^{d}$ be the vector of inputs and \\
$w_i:=\left(w_{i1},\cdots,w_{id}\right)^T\in{\mathbb R}^{d}$ be a parameter vector for the hidden unit $i$. The MLP function with $k$ hidden units can be written : 

\[
f_\theta(x)=\beta+\sum_{i=1}^k a_i\phi\left(w_i^Tx+b_i\right),
\]
with $\theta=\left(\beta,a_1,\cdots,a_k,b_1,\cdots,b_k,w_{11},\cdots,w_{1d},\cdots,w_{k1},\cdots,w_{kd}\right)$ the parameter vector of the model and  $\phi$ a bounded transfer function, usually a sigmo\"idal function. Note that we consider only real functions, extension to vectorial functions is straightforward but not discussed in this paper. 
Let $\Theta_k\subset \mathbb R^{k\times (d+2)+1}$ be a compact (i.e. closed and bounded) set of possible parameters, we consider regression model \({\mathcal S}=\left\{f_\theta(y,x),\ \theta\in\Theta_k\right\}\) with
\begin{equation}\label{model}
Y=f_{\theta}(X)+\varepsilon
\end{equation}
\(X\) is a random input variable and $\varepsilon$ is the noise of the model. Let $n$ be a strictly positive integer, we assume that the observed data \(\left(x_1,y_1\right),\cdots,\left(x_n,y_n\right)\) come from a true model \(\left(X_i,Y_i\right)_{i\in{\mathbb N},i>0}\) of which the true regression function is $f_{\theta^0}$, for an $\theta^0$ (possibly not unique) in the interior of $\Theta_k$. 
\subsection{Estimation of MLP regression model}
The main goal of non-linear regression is to give an estimation of the true parameter $\theta^0$ based on observations $\left((x_1,y_1),\cdots,(x_n,y_n)\right)$. This can be done by minimizing the Mean Square Error (MSE) function:
\begin{equation}\label{MSE}
E_n(\theta):=\frac{1}{n}\sum_{t=1}^{n}\left(y_t-f_{\theta}(x_t)\right)^2
\end{equation}
with respect to parameter vector $\theta\in\Theta_k$. The parameter vectors ${\hat \theta}_n$ realizing the minimum will be called Least Square Estimator (LSE). Note that parameters realizing the true distribution function may belong to a non-null dimension sub-manifold if the number of hidden units is overestimated. Suppose, for example, we have a multilayer perceptron with two hidden units and the true function $f_{\theta^0}$ is given by a perceptron with only one hidden unit, say $f_{\theta^0}=a_0\tanh(w_0x)$, with $x\in\mathbb R$. Then, any parameter $\theta$ in the set:
\[
\left\{\theta\left|w_2=w_1=w_0,b_2=b_1=0,a_1+a_2=a_0\right.\right\}
\]
realizes the function $f_{\theta^0}$. Hence, classical statistical theory for studying the LSE can not be applied because it requires the identification of the parameters (up to a permutation).

In the next section, we will compare MSE of over-determined models against MSE of the true model : 
\begin{equation}
\frac{1}{n}\sum_{t=1}^{n}\left(y_t-f_{\theta}(x_t)\right)^2-\frac{1}{n}\sum_{t=1}^{n}\left(y_t-f_{\theta^0}(x_t)\right)^2=E_n(\theta)-E_n(\theta^0).
\end{equation}

\section{A general bound for the MSE}
For an square integrable function $g(X,Y)$ the $L_2$  norm is: 
\begin{equation}\label{L2}
\Vert g(X,Y)\Vert_2:=\sqrt{\int g^2(x,y)dP(x,y)}.
\end{equation}
Now, for $\lambda>0$, let us define the generalized derivative function : 
\begin{equation}\label{gender}
d^\lambda_\theta(X,Y)=\frac{\frac{e^{-\lambda(Y-f_\theta(X))^2}-e^{-\lambda(Y-f_{\theta^0}(X))^2}}{e^{-\lambda(Y-f_{\theta^0}(X))^2}}}{\Vert\frac{e^{-\lambda(Y-f_\theta(X))^2}-e^{-\lambda(Y-f_{\theta^0}(X))^2}}{e^{-\lambda(Y-f_{\theta^0}(X))^2}}\Vert_2}=\frac{e^{-\lambda\left((Y-f_\theta(X))^2-(Y-f_{\theta^0}(X))^2\right)}-1}{\Vert e^{-\lambda\left((Y-f_\theta(X))^2-(Y-f_{\theta^0}(X))^2\right)}-1\Vert_2}
\end{equation}
and let us define \( \left(d^\lambda_\theta\right)_{-}(x,y)=\min\left\{0,d^\lambda_\theta(x,y)\right\}\).
Note that the generalized derivative function converges toward the derivative function if $\theta$ converges toward $\theta_0$. 

For now, let us assume that $d^\lambda_\theta$ is well defined, this point will be discuss later.  We can state the following inequality:

\textbf{\emph{Inequality}}:\\
{\it for} $\lambda>0$,
\begin{equation}
\sup_{\theta\in\Theta_k}n\cdot\left(E_n(\theta^0)-E_n(\theta)\right)\leq \frac{1}{2\lambda}\sup_{\theta\in\Theta_k}\frac{\sum_{i=1}^nd^\lambda_\theta(x_i,y_i)}{\sum_{i=1}^n\left(d^\lambda_\theta\right)^2_{-}(x_i,y_i)}
\end{equation}
\textbf{\emph{Proof}}:\\
We have
\[
\begin{array}{l}
n\cdot\left(E_n(\theta^0)-E_n(\theta)\right)=\\
\frac{1}{\lambda}\sum_{i=1}^n \log\left(1+\Vert\frac{e^{-\lambda(Y-f_\theta(X))^2}-e^{-\lambda(Y-f_{\theta^0}(X))^2}}{e^{-\lambda(Y-f_{\theta^0}(X))^2}}\Vert_2d^\lambda_\theta(x_i,y_i)\right)\\
\leq \sup_{0\leq p\leq \Vert\frac{e^{-\lambda(Y-f_\theta(X))^2}-e^{-\lambda(Y-f_{\theta^0}(X))^2}}{e^{-\lambda(Y-f_{\theta^0}(X))^2}}\Vert_2}\frac{1}{\lambda}\sum_{i=1}^n\log\left(1+pd^\lambda_\theta(x_i,y_i)\right)\\
\leq \sup_{p\geq 0}\frac{1}{\lambda}\left(p\sum_{i=1}^nd^\lambda_\theta(x_i,y_i)-\frac{p^2}{2}\sum_{i=1}^n\left(d^\lambda_\theta\right)^2_{-}(x_i,y_i)\right).
\end{array}
\]
Since for any real number $u$, \(\log(1+u)\leq u-\frac{1}{2}u^2_{-}\). Finally, replacing $p$ by the optimal value, we found
\[
\begin{array}{l}
n\cdot\left(E_n(\theta^0)-E_n(\theta)\right)\leq\frac{1}{2\lambda}\frac{\sum_{i=1}^nd^\lambda_\theta(x_i,y_i)}{\sum_{i=1}^n\left(d^\lambda_\theta\right)^2_{-}(x_i,y_i)}\\
\blacksquare
\end{array}
\]

This inequality allows to prove the tightness of $n\cdot\left(E_n(\theta^0)-E_n(\theta)\right)$ under simple assumptions. It is used in the next section to prove consistency of an estimator of the number of hidden units using penalized MSE criterion.
\section{Estimation of the number of hidden units.}
Let $k^0$ be the minimal number of hidden units needed to realize the true regression function $f_{\theta^0}$.  In this section, the set $\Theta$ of possible parameters will be set to
\[
\Theta=\cup_{k=1}^K\Theta_k,
\]
where $K$ is a, possibly huge, fixed constant: The maximum number of hidden units for MLP models. We define the minimum-penalized MSE estimator of $k^0$, as the minimizer $\hat k$ of 
\begin{equation}\label{PMSE}
T_n(k)=\min_{\theta\in\Theta}\left(E_n(\theta)+a_n(k)\right)
\end{equation}
Let us assume the following assumptions:
\begin{description}
\item{\bf (A1)} $a_n(.)$ is increasing, $n\cdot(a_n(k_1)-a_n(k_2))$ tends to infinity as $n$ tends to infinity, for any $k_1>k_2$ and $a_n(k)$ tends to 0 as $n$ tends to infinity for any $k$.
\item{\bf (A2)} It exists $\lambda>0$ so that \(\left\{d^\lambda_\theta,\theta\in\Theta\right\}\) is a Donsker class (see van der Vaart \cite{Vaart}) and appendix.
\end{description}
We now have:\\
\textbf{\emph{Theorem}}:\\
{\it Under {\bf (A1)} and {\bf (A2)}, $\hat k$ converges in probability to the true number of hidden units $k^0$.}\\
\textbf{\emph{Proof}}:\\
By applying the inequality,
\[
\begin{array}{l}
P(\hat k>k^0)\leq \sum_{k=k^0+1}^{K}P\left(T_n(k)\geq T_n(k^0)\right)=\\
\sum_{k=k^0+1}^{K}P\left(n\left(\sup_{\theta\in\Theta_{k^0}}E_n(\theta)-\sup_{\theta\in\Theta_k}E_n(\theta)\right)\geq n\left(a_n(k)-a_n(k^0)\right)\right)\leq\\
\sum_{k=k^0+1}^{K}P\left(\frac{1}{\lambda}\sup_{\theta\in\Theta_k}\frac{\sum_{i=1}^nd^\lambda_\theta(x_i,y_i)}{\sum_{i=1}^n\left(d^\lambda_\theta\right)^2_{-}(x_i,y_i)}\geq n\left(a_n(k)-a_n(k^0)\right)\right)
\end{array}
\]
Now, under {\bf (A2)} 
\[
sup_{\theta\in\Theta_k}\frac{1}{n}\left(\sum_{i=1}^nd^\lambda_\theta(x_i,y_i)\right)^2=O_P(1)
\]
where, $O_p(1)$ means bounded in probability. Moreover,  under {\bf (A2)}  the set $\left\{\left(d^\lambda_\theta(x_i,y_i)\right)^2\right\}$ is Glivenko-Cantelli (the set admits an uniform law of large numbers). Hence
\[
\inf_{\theta\in\Theta_k}\frac{1}{n}\sum_{i=1}^n\left(d^\lambda_\theta(x_i,y_i)\right)^2_{-}\stackrel{n\rightarrow\infty}{\longrightarrow}\inf_{\theta\in\Theta_k}\Vert\left(d^\lambda_\theta(X,Y)\right)_{-}\Vert_2^2
\]
But \(\inf_{\theta\in\Theta_k}\Vert\left(d^\lambda_\theta(X,Y)\right)_{-}\Vert_2>0\), since the random variable $d^\lambda_\theta(X,Y)$ is centered and $\Vert d^\lambda_\theta(X,Y)\Vert_2=1$.
Then, we get :
\[
\frac{1}{\lambda}\sup_{\theta\in\Theta_k}\frac{\sum_{i=1}^nd^\lambda_\theta(x_i,y_i)}{\sum_{i=1}^n\left(d^\lambda_\theta\right)^2_{-}(x_i,y_i)}=O_P(1)
\] 
and $P(\hat k>k^0)$ tends to 0 as $n$ tends to infinity. Finally,
\[
P(\hat k<k^0)\leq\sum_{k=1}^{k^0-1}P\left(\sup_{\theta\in\Theta_k}\frac{E_n(\theta)-E_n(\theta^0)}{n}\geq\frac{a_n(k)-a_n(k^0)}{n}\right)
\]
and \(\sup_{\theta\in\Theta_k}\frac{E_n(\theta)-E_n(\theta^0)}{n}\) converges in probability to 
\[
\sup_{\theta\in\Theta_k}E\left(E_n(\theta)-E_n(\theta^0)\right)<0
\] since $k<k^0$, so \(\hat k\stackrel{P}{\longrightarrow}k^0\) 
  $\blacksquare$

The assumption {\bf (A1)} is fairly standard for model selection, in the Gaussian case {\bf (A1)} will be fulfilled by the BIC criterion. The assumption  {\bf (A2)} is more difficult to check. First we note: 
\[
\begin{array}{l}
\left(e^{-\lambda\left((Y-f_\theta(X))^2-(Y-f_{\theta^0}(X))^2\right)}-1\right)^2=\\
e^{-2\lambda\left((Y-f_\theta(X))^2-(Y-f_{\theta^0}(X))^2\right)}-2e^{-\lambda\left((Y-f_\theta(X))^2-(Y-f_{\theta^0}(X))^2\right)}+1
\end{array}
\]
So, $d^\lambda_\theta$ is well defined if \(E\left[e^{-2\lambda\left((Y-f_\theta(X))^2-(Y-f_{\theta^0}(X))^2\right)}\right]<\infty\), but
\[
\begin{array}{l}
(Y-f_\theta(X))^2-(Y-f_{\theta^0}(X))^2=\\
(Y-f_{\theta^0}(X)+f_{\theta^0}(X)-f_\theta(X))^2-(Y-f_{\theta^0}(X))^2=\\
2\varepsilon(f_{\theta^0}(X)-f_\theta(X))+(f_{\theta^0}(X)-f_\theta(X))^2
\end{array}
\]
where $\varepsilon=Y-f_{\theta^0}(X)$ is the noise of the model. Since an MLP function is bounded,  $d^\lambda_\theta$ is well defined if $\lambda>0$ exists  such that $e^{\lambda\left|\varepsilon\right|}<\infty$ i.e. $\varepsilon$ admits exponential moments. Finally, using the same techniques of reparameterization as in Rynkiewicz \cite{Rynkiewicz}, assumption {\bf (A2)} can be shown to be true for MLP models with sigmo\"idal transfer functions, if the set of possible parameters $\Theta$ is compact.

\section{A little experiment}
The theoretical penalization terms of the previous section can be chosen among a wide range of functions (see condition {\bf A1}). In the sequel, a little experiment is conducted to assess the right rate of penalization to guess the ``true'' architecture of a model. 

Consider a simulated model:
\[
Z_t=F_{\theta^0}(X_t,Y_t)+\varepsilon_t, t=1,\cdots,n,
\]
with $\left((X_1,Y_1),\cdots,(X_n,Y_n)\right)$ i.i.d., $(X_{t},Y_{t})\sim{\cal N}\left(0_{\mathbb R^2},3\cdot I_2\right)$,\\
$\left(\varepsilon_1,\cdots,\varepsilon_n\right)$ i.i.d., $\varepsilon_t\sim{\cal U}\left[-1,1\right]$, the uniform law in $[-1;1]$ and
\begin{equation}\label{truemodel}
\begin{array}{l}
F_{\theta^0}(x,y)= \tanh(6\cdot x-2\cdot y)+2\cdot\tanh(8-x+3\cdot y)\\
-3\cdot tanh(2-6\cdot x-2\cdot y)+1.5.
\end{array}
\end{equation}
Here, the true model is an MLP with 2 inputs, 3 hidden units and one output. In order to avoid too long time of computation, the number of hidden units is assumed to be between $1$ and $10$.  We estimate the true architecture of the MLP according to (\ref{PMSE}).

First, let us write the  log-likelihood of the data as if the density of the noise would be Gaussian :
\begin{equation}
\begin{array}{l}
l_\theta\left(\left(\begin{array}{c}x_1\\y_1\\z_1\end{array}\right),\cdots,\left(\begin{array}{c}x_n\\y_n\\z_n\end{array}\right)\right)=-\frac{n}{2}\cdot\log(2\pi\sigma^2)\\
-\sum_{t=1}^n\frac{1}{2\sigma^2}\left(z_t-F_{\theta}(x_t,y_t)\right)^2+\sum_{t=1}^ng(x_t,y_t)
\end{array}
\end{equation}
Here, we assume that the variance of the noise $\sigma^2$ is a known constant and that the density of the explicative variables $(X,Y)$ is a function $g(.,.)$ independent of the parameter vector $\theta$.
Two classical criteria are:
\begin{equation}
\begin{array}{l}
AIC\,:\,\sum_{t=1}^n\frac{1}{\sigma^2}\left(z_t-F_{\theta}(x_t,y_t)\right)^2+2\cdot D+Cte\\
BIC\,:\,\sum_{t=1}^n\frac{1}{\sigma^2}\left(z_t-F_{\theta}(x_t,y_t)\right)^2+ D\cdot \log(n)+Cte\\
\end{array}
\end{equation}
where $D$ is the size of the parameter vector (the dimension of the model or the number of weights of the MLP) and $Cte$ a constant independent of the parameter $\theta$.

The optimization of the log-likelihood is done with respect to the parameter $\theta$, so maximimizing this quantity is equivalent to minimizing :
\begin{equation}
E_n(\theta)=\frac{1}{n}\cdot\sum_{t=1}^n\left(z_t-F_{\theta}(x_t,y_t)\right)^2
\end{equation}
Hence, an AIC like minimum-penalized MSE criterion would be:
\[
\frac{1}{n}\sum_{t=1}^n\left(z_t-F_{\theta}(x_t,y_t)\right)^2+\frac{2\sigma^2D}{n}
\]
and, a BIC like minimum-penalized MSE criterion would be:
\[
\frac{1}{n}\sum_{t=1}^n\left(z_t-F_{\theta}(x_t,y_t)\right)^2+ \frac{\sigma^2D\log(n)}{n}
\]
Note, that these criteria involve the knowledge of the variance of the noise $\sigma^2$. In a first experiment we will use the true variance of the noise, then this problem will be addressed in the following sections.

In the following sections, the optimization of MLP is done with the Broyden–Fletcher–Goldfarb–Shanno (BFGS) method. In order to avoid bad local minima, 10 random initializations of the weights are done for each estimation. 
\subsection{Model selection with $\sigma^2$ known}
We will compare 4 criteria, from the least penalized (AIC like) to the most penalized (Very Strong Penalization),
the following penalized criteria are assessed: 
\begin{itemize}
\item AIC like: $\frac{1}{n}\sum_{t=1}^n\left(z_t-F_{\theta}(x_t,y_t)\right)^2+\frac{2\sigma^2D}{n}$ 
\item BIC like: $\frac{1}{n}\sum_{t=1}^n\left(z_t-F_{\theta}(x_t,y_t)\right)^2+\frac{\sigma^2D\log(n)}{n}$
\item SP (Strong Penalization): $\frac{1}{n}\sum_{t=1}^n\left(z_t-F_{\theta}(x_t,y_t)\right)^2+\frac{\sigma^2D\sqrt{n}}{n}$
\item VSP (Very Strong Penalization): $\frac{1}{n}\sum_{t=1}^n\left(z_t-F_{\theta}(x_t,y_t)\right)^2+\frac{\sigma^2Dn^{3/4}}{n}$
\end{itemize}
We simulate $n=100$, $n=500$ and $n=1000$ data according to the true model (\ref{truemodel}), for each $n$ the experiment is repeated $100$ times.

The following architectures are chosen by the penalized criteria :
\begin{itemize}
\item n=100 
\[
\begin{array}{|c|c|c|c|c|c|c|c|c|c|c|c|}
\hline
&\mbox{nb h. units}&1&2&3&4&5&6&7&8&9&10\\
\hline
\mbox{AIC like}&\mbox{models sel.}&0&0&{\bf 47}&{\bf 35}&{\bf 10}&{\bf 2}&{\bf 3}&{\bf 3}&0&0\\
\hline
\mbox{BIC like}&\mbox{models sel.}&0&0&{\bf 100}&0&0&0&0&0&0&0\\
\hline
\mbox{SP}&\mbox{models sel.}&0&0&{\bf 100}&0&0&0&0&0&0&0\\
\hline
\mbox{VSP}&\mbox{models sel.}&0&{\bf 74}&{\bf 26}&0&0&0&0&0&0&0\\
\hline
\end{array}
\]

\item n=500 
\[
\begin{array}{|c|c|c|c|c|c|c|c|c|c|c|c|}
\hline
&\mbox{nb h. units}&1&2&3&4&5&6&7&8&9&10\\
\hline
\mbox{AIC like}&\mbox{models sel.}&0&0&{\bf 3}&{\bf 14}&{\bf 17}&{\bf 14}&{\bf 13}&{\bf 15}&{\bf 10}&{\bf 13}\\
\hline
\mbox{BIC like}&\mbox{models sel.}&0&0&{\bf 100}&0&0&0&0&0&0&0\\
\hline
\mbox{SP}&\mbox{models sel.}&0&0&{\bf 100}&0&0&0&0&0&0&0\\
\hline
\mbox{VSP}&\mbox{models sel.}&0&{\bf 1}&{\bf 99}&0&0&0&0&0&0&0\\
\hline
\end{array}
\]

\item n=1000 
\[
\begin{array}{|c|c|c|c|c|c|c|c|c|c|c|c|}
\hline
&\mbox{nb h. units}&1&2&3&4&5&6&7&8&9&10\\
\hline
\mbox{AIC like}&\mbox{models sel.}&0&0&0&0&0&0&0&{\bf 7}&{\bf 24}&{\bf 69}\\
\hline
\mbox{BIC like}&\mbox{models sel.}&0&0&{\bf 100}&0&0&0&0&0&0&0\\
\hline
\mbox{SP}&\mbox{models sel.}&0&0&{\bf 100}&0&0&0&0&0&0&0\\
\hline
\mbox{VSP}&\mbox{models sel.}&0&0&{\bf 100}&0&0&0&0&0&0&0\\
\hline
\end{array}
\]
\end{itemize}

The BIC like criterion and the Strong Penalization chose always the true architecture whatever the number of data. According to the theory, AIC like criterion is not consistent (see condition {\bf A1}) and the chosen architecture is always too large. The Very Strong penalization chose a too small architecture when the number of data is small ($n=100$), however it is a consistent criterion, so its behavior is correct for larger number of data ($n=500$ and $n=1000$).
This good results assume that the true variance of the noise is known, but for regression models this is never the case. A first idea may be to replace the unknown variance by the estimated one, the is done in the next section.

\subsection{Model selection using estimated variance ${\hat \sigma}^2$.}
The estimated variance ${\hat \sigma}^2$ is the mean square error of the model:
\[
{\hat \sigma}^2:=E_n(\hat \theta)=\frac{1}{n}\sum_{t=1}^n\left(z_t-F_{\hat \theta}(x_t,y_t)\right)^2
\] 
computed for the least square estimator $\hat \theta$. 
Hence, the comparison is done with the penalized criteria :
\begin{itemize}
\item AIC like: $\frac{1}{n}\sum_{t=1}^n\left(z_t-F_{\theta}(x_t,y_t)\right)^2+\frac{2{\hat \sigma}^2D}{n}$ 
\item BIC like: $\frac{1}{n}\sum_{t=1}^n\left(z_t-F_{\theta}(x_t,y_t)\right)^2+\frac{{\hat \sigma}^2D\log(n)}{n}$
\item SP (Strong Penalization): $\frac{1}{n}\sum_{t=1}^n\left(z_t-F_{\theta}(x_t,y_t)\right)^2+\frac{{\hat \sigma}^2D\sqrt{n}}{n}$
\item VSP (Very Strong Penalization): $\frac{1}{n}\sum_{t=1}^n\left(z_t-F_{\theta}(x_t,y_t)\right)^2+\frac{{\hat \sigma}^2Dn^{3/4}}{n}$
\end{itemize}

The following architectures are chosen by the penalized criteria: 
\begin{itemize}
\item n=100 
\[
\begin{array}{|c|c|c|c|c|c|c|c|c|c|c|c|}
\hline
&\mbox{nb h. units}&1&2&3&4&5&6&7&8&9&10\\
\hline
\mbox{AIC like}&\mbox{models sel.}&0&0&0&0&0&{\bf 1}&0&{\bf 4}&{\bf 27}&{\bf 68}\\
\hline
\mbox{BIC like}&\mbox{models sel.}&0&0&{\bf 22}&{\bf 6}&{\bf 3}&{\bf 1}&{\bf 3}&{\bf 5}&{\bf 15}&{\bf 45}\\
\hline
\mbox{SP}&\mbox{models sel.}&0&0&{\bf 85}&{\bf 3}&{\bf 1}&0&0&0&{\bf 1}&{\bf 10}\\
\hline
\mbox{VSP}&\mbox{models sel.}&0&0&{\bf 100}&0&0&0&0&0&0&0\\
\hline
\end{array}
\]

\item n=500 
\[
\begin{array}{|c|c|c|c|c|c|c|c|c|c|c|c|}
\hline
&\mbox{nb h. units}&1&2&3&4&5&6&7&8&9&10\\
\hline
\mbox{AIC like}&\mbox{models sel.}&0&0&0&{\bf 4}&{\bf 4}&{\bf 6}&{\bf 9}&{\bf 23}&{\bf 19}&{\bf 35}\\
\hline
\mbox{BIC like}&\mbox{models sel.}&0&0&{\bf 100}&0&0&0&0&0&0&0\\
\hline
\mbox{SP}&\mbox{models sel.}&0&0&{\bf 100}&0&0&0&0&0&0&0\\
\hline
\mbox{VSP}&\mbox{models sel.}&0&0&{\bf 100}&0&0&0&0&0&0&0\\
\hline
\end{array}
\]

\item n=1000 
\[
\begin{array}{|c|c|c|c|c|c|c|c|c|c|c|c|}
\hline
&\mbox{nb h. units}&1&2&3&4&5&6&7&8&9&10\\
\hline
\mbox{AIC like}&\mbox{models sel.}&0&0&0&0&{\bf 3}&{\bf 7}&{\bf 12}&{\bf 16}&{\bf 29}&{\bf 33}\\
\hline
\mbox{BIC like}&\mbox{models sel.}&0&0&{\bf 100}&0&0&0&0&0&0&0\\
\hline
\mbox{SP}&\mbox{models sel.}&0&0&{\bf 100}&0&0&0&0&0&0&0\\
\hline
\mbox{VSP}&\mbox{models sel.}&0&0&{\bf 100}&0&0&0&0&0&0&0\\
\hline
\end{array}
\]
\end{itemize}
As usual the AIC like criterion misbehaves like in the previous section.
But, for a small number of data ($n=100$) the use of an estimation of the variance of the noise instead of the true one leads to overestimation of the number of hidden unit for the BIC like criterion and the strong penalization.
The explanation is that the variance of the noise is underestimated for large number of hidden units and so the penalized criterion.
This drawback disappears for larger number of data ($n=500$ and $n=1000$) because the estimation of the variance becomes better.  Despite that the Very Strong Penalized criterion seems to guess the true architecture whatever the number of data, to plug estimated variance instead of the true one seems maybe to be a too naive approach for small number of data.  As the goal of penalized criterion is to compare models, we could use the logarithm of the mean square error instead of the mean square error itself, hence the lack of the true variance is no more a problem because this number is simplified in the difference of the logarithms. This approach is studied in the next section.

\subsection{Model selection using logarithm of mean square error}

Choosing between two number of hidden units, $k_1$ and $k_2$ is the result of a comparison between $\min_{\theta_1\in\Theta_{k_1}}\left(E_n(\theta_1)+a_n(k_1)\right)$ on one hand and $\min_{\theta_2\in\Theta_{k_2}}\left(E_n(\theta)+a_n(k_2)\right)$ on the other hand, where $a_n(k)$ is the penalization term. So, the results of the comparison may be changed if we consider $C\cdot E_n(\theta)$ instead of $E_n(\theta)$ when $C$ is a very big (or very small) constant. But, if we compare $\min_{\theta_1\in\Theta_{k_1}}\left(\log\left(E_n(\theta_1)\right)+a_n(k_1)\right)$ and\\
$\min_{\theta_2\in\Theta_{k_2}}\left(\log\left(E_n(\theta_2)\right)+a_n(k_2)\right)$ the results of the comparison is the same if we change $E_n(\theta)$ in  $C\cdot E_n(\theta)$.

Moreover,  $E_n(\theta)$ may be seen as an approximation of the variance of the noise $\sigma^2$, so $E_n(\theta)=\sigma^2\cdot(1-\varepsilon(\theta))$ and 
\[
\log\left(E_n(\theta)\right)\simeq\log(\sigma^2)-\varepsilon(\theta)+o(\varepsilon(\theta))
\]

The term $\varepsilon(\theta)$ may be seen as the normalized term of overfitting of the model. Finally we get:
\[
\begin{array}{l}
\min_{\theta\in\Theta}\left(\log\left(E_n(\theta)\right)+a_n(k_1)\right)-\min_{\theta\in\Theta}\left(\log\left(E_n(\theta)\right)+a_n(k_2)\right)\simeq\\
\varepsilon(\theta_2)-\varepsilon(\theta_1)+a_n(k_2)-a_n(k_1)
\end{array}
\]  
and the penalization term plays fully is role of compensation of the ``normalized'' overfitting.

The results on our little experiment with these new criteria are the following:
\begin{itemize}
\item n=100 
\[
\begin{array}{|c|c|c|c|c|c|c|c|c|c|c|c|}
\hline
&\mbox{nb h. units}&1&2&3&4&5&6&7&8&9&10\\
\hline
\mbox{AIC like}&\mbox{models sel.}&0&0&{\bf 3}&{\bf 7}&{\bf 4}&{\bf 11}&{\bf 9}&{\bf 13}&{\bf 24}&{\bf 29}\\
\hline
\mbox{BIC like}&\mbox{models sel.}&0&0&{\bf 96}&{\bf 4}&0&0&0&0&0&0\\
\hline
\mbox{SP}&\mbox{models sel.}&0&0&{\bf 100}&0&0&0&0&0&0&0\\
\hline
\mbox{VSP}&\mbox{models sel.}&{\bf 73}&{\bf 27}&0&0&0&0&0&0&0&0\\
\hline
\end{array}
\]

\item n=500 
\[
\begin{array}{|c|c|c|c|c|c|c|c|c|c|c|c|}
\hline
&\mbox{nb h. units}&1&2&3&4&5&6&7&8&9&10\\
\hline
\mbox{AIC like}&\mbox{models sel.}&0&0&{\bf 3}&{\bf 6}&{\bf 8}&{\bf 13}&{\bf 18}&{\bf 10}&{\bf 17}&{\bf 25}\\
\hline
\mbox{BIC like}&\mbox{models sel.}&0&0&{\bf 99}&{\bf 1}&0&0&0&0&0&0\\
\hline
\mbox{SP}&\mbox{models sel.}&0&0&{\bf 100}&0&0&0&0&0&0&0\\
\hline
\mbox{VSP}&\mbox{models sel.}&0&{\bf 69}&{\bf 31}&0&0&0&0&0&0&0\\
\hline
\end{array}
\]

\item n=1000 
\[
\begin{array}{|c|c|c|c|c|c|c|c|c|c|c|c|}
\hline
&\mbox{nb h. units}&1&2&3&4&5&6&7&8&9&10\\
\hline
\mbox{AIC like}&\mbox{models sel.}&0&0&{\bf 1}&{\bf 2}&{\bf 7}&{\bf 13}&{\bf 16}&{\bf 21}&{\bf 19}&{\bf 21}\\
\hline
\mbox{BIC like}&\mbox{models sel.}&0&0&{\bf 100}&0&0&0&0&0&0&0\\
\hline
\mbox{SP}&\mbox{models sel.}&0&0&{\bf 100}&0&0&0&0&0&0&0\\
\hline
\mbox{VSP}&\mbox{models sel.}&0&{\bf 2}&{\bf 98}&0&0&0&0&0&0&0\\
\hline
\end{array}
\]
\end{itemize} 

We can see that this method yields very good results, whatever the number of data, for BIC like criterion and the Strong Penalization without knowing the true variance of the noise. Strong Penalization seems even to be a little better than BIC like criterion.

\section{Conclusion}
MLP regression is widely used and always a very competitive method (see Osowski et al. \cite{Osowski}), but theoretical justification is lacking for determining the true architecture and especially the number of hidden units.  Indeed, the classical asymptotic theory fails when the model is not identifiable. In this paper, we prove an inequality showing that overfitting of MLP is moderate if the noise admits exponential moments and the parameters of the model are a priori bounded. This bound justifies  the use of penalized criteria in order to fit the architecture of MLP models in the framework of regression without knowing the density of the noise.
Hence, The user can select the true number of hidden units thanks to  penalized criteria, of the form
\[
\begin{array}{l}
E_n(\theta)+a_n(k)\\
\mbox{or}\\
\ln(E_n(\theta))+a_n(k)
\end{array}
\]

If the penalization term $a_n(k)$ is well calibrated ($\frac{C\cdot k}{n}<a_n(k)<C\cdot k$), the true number of hidden units  will be automatically selected if  $n$ is large enough.
A little experiment suggests that a good choice of penalization seems to be the middle of the possible range: $a_n(k)=\frac{C\cdot k}{\sqrt{n}}$. The use of the logarithm of the mean square error $E_n(\theta)$ is an easy way to avoid to know the true variance of the noise. A further question could be to know if this empirical finding for the tuning of the penalization term can be justified theoretically.

Note that, this paper was only concerned with the identification of the true model. The point is more to have an idea of the complexity of the model determining the data than to have a predictive model. However, if there are enough data, the true model will also be the best predictor. Hence, the spirit of the studied criterion is then very different from the approach used in Extreme Learning Machine (ELM) which provides a good predictive model at extremely fast learning speed. Indeed, in ELM, all the hidden node parameters are independent from the target functions or the training datasets so the weights between input and hidden node are no more considered as parameters. Such method is really not concerned by model identification but only by predictive power and speed of computing, for example it is possible to use ELM even if the relation between inputs and ouput is linear. Finally, If ELM seems to work well in pratice, the theoritical justification of the superiority of such method is still lacking. 

Another issue may be when the unknow regression function is not represented by an MLP at all, in such case there is no more over-determination and the asymptotic behaviour of the model is more underfitting than overfitting. The overfitting occures only when we consider finite number of data and the theory to deal with such problem is the very difficult non-asymptotic statistical theory as in Massart \cite{Massart}. As far as we know this theory gives no hints for choosing the number of hidden units for finite number of data.

\section*{Appendix}
A Donsker class is a notion from the ``empirical processes theory''. This theory deals with ``law of large number'' and ``asymptotic normality'' for set of functions.  Basically, a Donsker class is a set of functions which is ``not too big''.

Let $X_1,\cdots,X_n$ be a random sample from a probability distribution $P$. The empirical distribution is the discrete uniform measure on the observations. We denote it by ${\mathbb P}_n=\frac{1}{n}\sum_{i=1}^n \delta_{X_i}$, where $\delta_x$ is the probability distribution that is degenerate at $x$. Given a function $f$, we write ${\mathbb P}_n f$ for the expectation of $f$ under the empirical measure and $Pf$ for the expectation under $P$. Thus
\[
{\mathbb P}_n f=\frac{1}{n}\sum_{i=1}^n f(X_i)\mbox{ and }Pf=\int fdP.
\]
The empirical process evaluated at $f$ is defined as ${\mathbb G}_nf=\sqrt{n}\left({\mathbb P}_n f-Pf\right)$.
Consider a set of functions $\mathcal{F}$  endowed with the $L_2$ norm $\left\Vert \cdot\right\Vert$ (see (\ref{L2})). For every $\varepsilon>0$, we define an $\varepsilon$-bracket by \\
 $\left[l,\, u\right]=\left\{ f\in\mathcal{F},\, l\leq f\leq u\right\} $ such that $\left\Vert u-l\right\Vert <\varepsilon$. The $\varepsilon$-bracketing entropy is
\[\mathcal{H}_{\left[\cdot\right]}\left(\varepsilon,\mathcal{F},\left\Vert \cdot\right\Vert\right)=\ln \left(\mathcal{N}_{\left[\cdot\right]}\left(\varepsilon,\mathcal{F},\left\Vert \cdot\right\Vert \right)\right),\]
where $\mathcal{N}_{\left[\cdot\right]}\left(\varepsilon,\mathcal{F},\left\Vert \cdot\right\Vert \right)$
is the minimum number of $\varepsilon$-brackets necessary to cover $\mathcal{F}$. $\mathcal{N}_{\left[\cdot\right]}\left(\varepsilon,\mathcal{F},\left\Vert \cdot\right\Vert \right)$ is also called ``covering number''.

The class $\cal F$ of functions is called Donsker if  the covering number, which depends of the diameter $\varepsilon$ of the balls, is smaller than order $e^{\frac{1}{\varepsilon^2}}$ when $\varepsilon$ goes to $0$. Then, the sequence of processes $\left\{{\mathbb G}_nf,f\in \mathcal{F}\right\}$ converges in distribution to a tight limit process. Finally, if ${\cal F}=\left\{d^\lambda_\theta,\theta\in\Theta\right\}$ is a Donsker class, the sequence of processes 
\[
\left(\frac{1}{\sqrt{n}}\left(\sum_{i=1}^nd^\lambda_\theta(X_i,Y_i)\right)\right)_{n\in\mathbb N}
\]
converges in distribution to a tight random Gaussien process.
We get then the key property needed for the demonstration of the theorem:
\[
sup_{\theta\in\Theta_k}\frac{1}{n}\left(\sum_{i=1}^nd^\lambda_\theta(x_i,y_i)\right)^2=O_P(1)
\]

\begin{footnotesize}

\end{footnotesize}
\end{document}